\documentclass{article}
\usepackage{RR}

\usepackage{amsmath,amsfonts,amssymb,latexsym,stmaryrd}
\usepackage[latin1]{inputenc}
\usepackage[LGR,T1]{fontenc}
\usepackage[plain]{algorithm}
\usepackage{algorithmic}
\usepackage{url}

\newcommand{\R}      {\mathbb R}

\newcommand{\ot}        {\leftarrow}
\newcommand{\carre}     {\hfill$\Box$}

\newcommand{\sfX}  {{\mathsf X}}

\newcommand{\sfs}  {{\mathsf s}}

\newcommand{\sfv}  {{\mathsf v}}
\newcommand{\sfw}  {{\mathsf w}}
\newcommand{\sfx}  {{\mathsf x}}
\newcommand{\sfy}  {{\mathsf y}}

\newcommand{\bfa}  {{\mathbf a}}
\newcommand{\bfb}  {{\mathbf b}}

\newcommand{\BB}   {{\cal B}}

\newcommand{\NN}   {{\cal N}}

\newcommand{\UU}   {{\cal U}}

\newcommand{\rmd}   {{{\textrm{\upshape d}}}}

\newcommand{\indic}{{\rm\bf1}}
\renewcommand{\epsilon}{\varepsilon}

\def\dobm{
    \copy1\kern-\wd1\kern0.05ex\copy1\kern-\wd1\kern0.05ex\box1}

\newcommand{\fenumi}{\rm({\it i}\/)}
\newcommand{\fenumii}{\rm({\it ii}\/)}
\newcommand{\fenumiii}{\rm({\it iii}\/)}

\newcommand{\eqdef}     {\stackrel{{\textrm{\rm\tiny def}}}{=}}
\newtheorem{theorem}      {Theorem}[section]
\newtheorem{theorem*}     {theorem}
\newtheorem{proposition}  [theorem]{Proposition}

\newtheorem{lemma}        [theorem]{Lemma}

\newcommand{\proof}        {\paragraph{Proof}}

\newcommand{\law}       {{\rm law}}

\newcommand{\piprop}  {\pi^{\textrm{\rm\tiny prop}}}
\newcommand{\pipropl}  {\pi^{\ell,\textrm{\rm\tiny prop}}}

\newcommand{\bbra}[1]{\llbracket #1 \rrbracket}

\begin{document}

\RRdate{October 2006}
\RRNo{xxxx}
\RRauthor{
Fabien Campillo 
\thanks{INRIA/IRISA, Rennes, \protect\url{Fabien.Campillo@inria.fr}}  
and
Vivien Rossi 
\thanks{IURC, University of Montpellier I -- \protect\url{Viven.Rossi@iurc.montp.inserm.fr}}
\thanks{The research of the second author was done during 
                  a postdoctoral stay at the INRIA/IRISA, Rennes.}}
\authorhead{F. Campillo \&\ V. Rossi} 
\RRetitle{Parallel and interacting\\ Markov chains Monte Carlo method} 
\titlehead{Parallel and interacting MCMC's} 
\RRtitle{Méthode de Monte Carlo par chaînes de Markov\\ en parallèle et en interaction} 
\RRabstract{In many situations it is important to be able to propose $N$ independent realizations of a given distribution law. We propose a strategy for making $N$ parallel Monte Carlo Markov Chains (MCMC) interact in order to get an approximation of an independent $N$-sample of a given target law. In this method each individual chain proposes candidates for all other chains. We prove that the set of interacting chains is itself a MCMC method for the product of $N$ target measures. Compared to independent parallel chains this method is more time consuming, but we show through concrete examples that it possesses many advantages: it can speed up convergence toward the target law as well as handle the multi-modal case.}

\RRresume{Dans de nombreuses situations il est important de pouvoir disposer de $N$ réalisations indépendantes d'une loi donnée. Notre but est de développer une stratégie d'interaction de $N$ méthodes de Monte Carlo par Chaîne de Markov (MCCM) dans le but de proposer une approximation d'un  échantillon indépendant de taille $N$ d'une loi cible donnée. L'idée est que chaque chaîne propose un candidat pour elle-même mais également pour toutes les autres chaînes. On montre que l'ensemble de ces $N$ chaînes en interaction est lui-même une méthode MCCM pour le produit de $N$ mesures cibles. Cette approche est naturellement plus coûteuse que $N$ chaînes indépendantes, on montre toutefois au travers d'exemples concrets qu'elle possède plusieurs avantages\,: elle peut sensiblement accélérer la convergence vers la loi cible, elle permet également d'appréhender le cas multimodal.
}
\RRkeyword{Markov chain Monte Carlo method, Metropolis-Hastings, interacting chains, particle approximation}
\RRmotcle{méthode de Monte Carlo par chaîne de Markov, Metropolis-Hastings, chaînes en interaction, approximation particulaire}

\RRprojets{Aspi} 
\RRtheme{\THNum} 

\URRennes

\makeRR 


\cleardoublepage

\tableofcontents

\cleardoublepage


\section{Introduction}

Markov chain Monte Carlo (MCMC) algorithms \cite{tierney1994a,gilks1995a,robert1996a} allows us  to draw  samples from a probability distribution $\pi(x)\,\rmd x$ known up to a multiplicative constant. They consist in sequentially simulating a single Markov chain whose limit distribution is $\pi(x)\,\rmd x$. There exist many techniques to speed up the convergence toward the target distribution by improving the mixing properties of the chain \cite{gilks1995b}. Moreover, special attention should be given to the convergence diagnosis of this method \cite{brooks1998a,cowles1996a,kass1998a}.

\smallskip

An alternative is to run many  Markov chains in parallel. The simplest multiple chain algorithm is to make use of parallel independent chains \cite{gelman1992a}. The recommendations concerning this idea seem contradictory in the literature (cf. the ``many short runs'' \textit{vs} ``one long run'' debate described in \cite{geyer1991a}). We can note with \cite{geyer1992a} and \cite[\S\,6.5]{robert1996a} that independent parallel chains could be a poor idea: among these chains some may not converge, so one long chain could be preferable to many short ones. Moreover, many parallel independent chains can artificially exhibit a more robust behavior which does not correspond to a real convergence of the algorithm.

\smallskip

In practice one however make use of several chains in parallel. It is then tempting to exchange information between these chains to improve mixing properties of the MCMC samplers \cite{chauveau2001a,chauveau2002a, laskey2003a, chao2003a, drugan2004a, drugan2005a}. A general framework of ``Population Monte Carlo'' has been proposed in this context \cite{iba2001a,mengersen2003a,cappe2004a}. 
In this paper we propose an interacting method between parallel chains which provides an independent sample from the target distribution. Contrary to papers previously cited, the proposal law n our work is given and does not adapt itself to the previous simulations. Hence, the problem of the choice of this law still remains.

\smallskip

The Metropolis-Hastings (MH) algorithm and its theoretical properties are presented in section \ref{sec.MH}. The corresponding Metropolis within Gibbs (MwG) algorithm and its theoretical properties are presented in section \ref{sec.para.int.MwG}. In Section \ref{sec.tests}, two simple numerical examples illustrate how the introduction of interactions can speed up the convergence and handle multi-modal cases.

\section[Parallel/interacting MH algorithm]
        {Parallel/interacting Metropolis Hastings (MH) algorithm}
\label{sec.MH}

Consider a target density law $\pi(x)$ defined on $(\R^n,\BB(\R^n))$ and a proposal kernel density $\piprop(y|x)$. We propose  a method for sampling $N$ independent values $X^1,\dots,X^N\in\R^n$ of the law  $\pi(x)\,\rmd x$. 

\paragraph{Notations:} {\it Let
\[
 X = X^{1:N}
   = X_{1:n} \in\R^{n\times N} \,,
\]
so that $X_\ell\in\R^N$ and $X^i\in\R^n$ (the same for $Y$ and $Z$); $x \in\R^n$ so that $x_\ell\in\R$ (the same for $y$ and $z$); $\xi,\xi'\in\R$.
Here  $X^{1:N} = (X^1,\dots,X^N)$ and  $X_{1:n} = (X_1,\dots,X_n)$. We also define $\neg\ell = \{1,\dots,n\} \setminus \{\ell\}$. Note that the structure of the matrix $X$ is:
\[
  X
  =
  \begin{array}{cc}
    \begin{array}{ccccc}
    \hphantom{X_1^1}  & \hphantom{\cdots} & X^i  & \hphantom{\cdots} & \hphantom{X^N_1}
   \\
    \hphantom{X_1^1}  & \hphantom{\cdots} & \ \uparrow \   & \hphantom{\cdots} & \hphantom{X^N_1}
  \end{array}
  \\
  \\
  \left[
  \begin{array}{ccccc}
    X_1^1  & \cdots & X^i_1  & \cdots & X^N_1  \\
    \vdots &        & \vdots &        & \vdots \\
    X^1_\ell & \cdots & X^i_\ell & \cdots & X^N_\ell    \\
    \vdots &        & \vdots &        & \vdots \\
    X^1_n & \cdots & X^i_n & \cdots & X^N_n
  \end{array}
  \right]
  &
  \begin{array}{c}
      \\
     \\
    \to\ X_\ell  \\
     \\
     \\
  \end{array}
  \\
  \vphantom{sdsdfdfsfsd}
  \\
  \vphantom{sdsdfdfsfsd}
  \\
  \vphantom{sdsdfdfsfsd}
  \end{array}
  \,.
\]
}

\subsection{The algorithm}
\label{sec.algo.mh}

We describe the Markov chain $\{\sfX^{(k)}\}_{k\geq 0}$ over $\R^{n\times N}$ corresponding the MH algorithm. It consists in $N$ \emph{mutually dependent} realizations $\sfX^{i,(k)}$ ($i=1,\dots,N$) of the state variable and its limit distribution will be
\[
 \Pi({\rmd} X) 
 \eqdef 
 \pi(X^1)\,\rmd X^1 \cdots \pi(X^N)\,\rmd X^N
 \,.
\]

\medskip

We detail an iteration $\sfX^{(k)}=X\to \sfX^{(k+1)}=Z$ of the MH algorithm. The $N$ vectors are updated sequentially: 
\[
  [X^{1:N}]
  \to
  [Z^1 X^{2:N}]
  \to
  [Z^{1:2} X^{3:N}]
  \cdots
  [Z^{1:N-1} X^{N}]
  \to
  [Z^{1:N}]
  \,.
\]
At sub-iteration ``$i$\,'', that is $[Z^{1:i-1} X^{i:N}] \to  [Z^{1:i} X^{i+1:N}]$, we simulate $Z^i$ in two steps:
\begin{description}

\item{\it Proposal step:}
independently one from the other, each chain $j=1\cdots N$ proposes a candidate $Y^j\in\R^n$ according to the proposal kernel starting from its current position, i.e.
\begin{align*}
   Y^j \sim \piprop_{i,j}(y|Z^{1:i-1},X^i,X^{i+1:N})\,\rmd y\,.
\end{align*}
Note that the candidates $Y^j$ depend also on $i$. We will use a lighter notation:
\begin{align}
\label{eq.proposal}
   \piprop_{i,j}(y|X^i)
   =
   \piprop_{i,j}(y|Z^{1:i-1},X^i,X^{i+1:N})
   \,.
\end{align}

\item{\it Selection step:}
We can chose among these $N$ candidates $Y^{1:N}$ or stay at $X^i$ according to the multinomial law:
\begin{align*}
   Z^i
   \ot
   \left\{\begin{array}{ll}
      Y^1 &\textrm{with probability } \frac{1}{N} \,\alpha^{i,1}(X^i,Y^1)\,,
      \\
      \ \vdots
      \\
      Y^{N} &\textrm{with probability } \frac{1}{N}\, \alpha^{i,N}(X^i,Y^{N})\,,
      \\[0.3em]
      X^i
      &\textrm{with probability } 
         \tilde\rho^i(X^i,Y)
   \end{array}\right.
\end{align*}
where the acceptance probabilities are  
\begin{align*}
   \alpha^{i,j}(x,y) 
   &\eqdef 
   \frac{\pi(y)}{\pi(x)}\,\frac{\piprop_{i,j}(x|y)}{\piprop_{i,j}(y|x)}
   \wedge 1
   \,,
\\
   \tilde\rho^i(X^i,Y)
   &\eqdef 
   1-\frac{1}{N} \sum_{j=1}^N \alpha^{i,j}(X^i,Y^j)
   \,.
\end{align*}
\end{description}
The final algorithm is depicted in Algorithm~\ref{algo.mh.1}.

\begin{algorithm}
\caption{\sl Parallel/interacting MH algorithm.}
\label{algo.mh.1}
\begin{center}
\begin{minipage}{12cm}
\hrulefill\\[-1em]
\mbox{}
\begin{algorithmic}
\STATE choose $X\in\R^{n\times N}$  
\FOR {$k=1,2,\dots$}
  \FOR {$i=1:N$}
    \FOR {$j=1:N$}
       \STATE $Y^j\sim \piprop_{i,j}(y|X^i)\,\rmd y$  
       \STATE $\alpha^j 
                \ot
                [ \pi(Y^j) \, \piprop_{i,j}(X^i|Y^j)]
                 / [\pi(X^i)\,\piprop_{i,j}(Y^j|X^i)]  \wedge 1$
    \ENDFOR
    \STATE $\tilde\rho \ot 1 - \frac{1}{N}\sum_{j=1}^N\alpha^j$
    \STATE $X^i
      \ot
      \left\{\begin{array}{ll}
        Y^1 &\textrm{with probability } \alpha^1/N
        \\[-0.4em]
        \ \vdots
        \\
        Y^{N} &\textrm{with probability } \alpha^N/N
        \\[0.2em]
        X^i
        &\textrm{with probability } 
           \tilde\rho
   \end{array}\right.$
  \ENDFOR
\ENDFOR
\end{algorithmic}
\hrulefill
\end{minipage}
\end{center}
\end{algorithm}

\subsection{Description of the MH kernel}

\begin{lemma}
\label{lemma.kernel.MH}
The Markov kernel associated with the MH procedure described in Section
\ref{sec.algo.mh} is
\begin{align}
\label{eq.noyau.mh.P}
   P(X;\rmd Z)
   \eqdef
   P^1(X^{1:N};\rmd Z^1) \; 
   P^2(Z^1,X^{2:N};\rmd Z^2) \cdots
   P^N(Z^{1:N-1},X^N;\rmd Z^N) 
\end{align}
where
\begin{align}
\label{eq.noyau.mh.P.i}
   P^i(Z^{1:i-1},X^{i:N};\rmd z) 
   &
   \eqdef
   \frac{1}{N}\sum_{j=1}^N \alpha^{i,j}(X^i,z) \; \piprop_{i,j}(z|X^i) \;\rmd z 
   +
   \rho^i(X^i) \; \delta_{X^i}(\rmd z)
   \,.
\end{align}
Acceptation probability is
\begin{align}
\label{eq.noyau.mh.alpha.ij}
   \alpha^{i,j}(x,z)
   &
   \eqdef
   \left\{\begin{array}{ll}
     r^{i,j}(x,z) \wedge 1 & \textrm{if }(x,z)\in R^{i,j} \,,\\
     0                     & \textrm{otherwise,}
   \end{array}\right.
\\[0.8em]
\label{eq.noyau.mh.r.ij}
   r^{i,j}(x,z)
   &\eqdef
   \frac{\pi(z)}{\pi(x)}\,\frac{\piprop_{i,j}(x|z)}{\piprop_{i,j}(z|x)}
   \,,
\\[0.8em]
\label{eq.noyau.mh.rho.i}
   \rho^i(x)
   &
   \eqdef
   1
   -
   \frac{1}{N}
   \sum_{j=1}^N \int_\R \alpha^{i,j}(x,z) \; \piprop_{i,j}(z|x) \;\rmd z 
   \,.
\end{align}
The set $R^{i,j}$ is defined by:
\begin{align*}
   R^{i,j}
   \eqdef
   \big\{
     (x,z)\in \R^n\times\R^n   \,;\,
    \pi(z) \; \piprop_{i,j}(x|z) > 0  
    \  \textrm{ and } \ 
    \pi(x) \; \piprop_{i,j}(z|x) > 0
    \big\} \,.
\end{align*}
Note that the functions  $\alpha^{i,j}(x,z)$, $\rho^i(x)$, $r^{i,j}(x,z)$ and the set $R^{i,j}$ depend on $Z^{1:i-1}$ and $X^{i:N}$.

The measures
\begin{align*}
  \nu(\rmd x\times\rmd z)
  &=
  \pi(z)\,\piprop_{i,j}(x|z) \,\rmd z\,\rmd x
  \,,
  &
  \nu^T(\rmd x\times\rmd z)
  &=
  \pi(x)\,\piprop_{i,j}(z|x) \,\rmd z\,\rmd x
\end{align*}
are mutually absolutely continuous over $R ^{i,j}$ and mutually singular on the complementary set $[R^{i,j}]^c$. The set $R ^{i,j}$ is unique, up to the $\nu$ and $\nu^T$ negligible sets, and symmetric, i.e. $(x,z)\in R^{i,j}\Rightarrow (z,x)\in R^{i,j}$.
\end{lemma}

\bigskip

\proof
This construction follows the general setup proposed by Luke Tierney in \cite{tierney1998a}.
We now derive the probability kernel associated with the iteration described in the previous subsection \ref{sec.algo.mh}. The kernel $P^i(Z^{1:i-1},X^{i:N};\rmd z)$ is the composition of a proposition kernel and of a selection kernel:
\begin{align*}
   &
   P^i(Z^{1:i-1},X^{i:N};\rmd z) 
   =
   \int_{Y^{1:N}}\!\!
     S^i(Z^{1:i-1},X^{i:N},Y^{1:N};\rmd z) \; 
     Q^i(Z^{1:i-1},X^{i:N};\rmd Y^{1:N})  
\end{align*}
which consists in proposing independently  $N$ candidates $Y^{1:N}$ sampled from the density proposition, i.e.
\begin{align*}
   Q^i(Z^{1:i-1},X^{i:N};\rmd Y^{1:N}) 
   &
   \eqdef
   \prod_{k=1}^N \piprop_{i,k}(Y^k|X^i)\,\rmd Y^k 
\end{align*}
then to select among these candidates or to stay at $X^i$ with the MH acceptance probability, i.e.
\begin{align*}
   S^i(Z^{1:i-1},X^{i:N},Y^{1:N};\rmd z)
   &
   \eqdef
   \frac{1}{N} \sum_{j=1}^N \alpha^{i,j}(X^i,Y^j)\,\delta_{Y^j}(\rmd z)
   +
   \tilde\rho^i(X^i,Y)\,\delta_{X^i}(\rmd z)
   \,.
\end{align*}
Hence:
\begin{align*}
   &
   P^i(Z^{1:i-1},X^{i:N};\rmd z) 
   =
   \\
   &
   \qquad\qquad
   =
   \frac{1}{N}\sum_{j=1}^N
   \int_{Y^{1:N}}
      \alpha^{i,j}(X^i,Y^j)\,\delta_{Y^j}(\rmd z)
      \;
      \Big\{\prod_{k=1}^N \piprop_{i,k}(Y^k|X^i)\,\rmd Y^k\Big\} 
  \\
  &\qquad\qquad\qquad
   +
   \int_{Y^{1:N}}
     \tilde\rho^i(X^i,Y)\,\delta_{X^i}(\rmd z)
     \;
     \Big\{\prod_{k=1}^N \piprop_{i,k}(Y^k|X^i)\,\rmd Y^k\Big\} 
  = A_1 + A_2
\end{align*}
and
\begin{align*}
   A_1
   &
   =
   \frac{1}{N}\sum_{j=1}^N
   \int_{Y^j}
      \alpha^{i,j}(X^i,Y^j)\,\delta_{Y^j}(\rmd z)\,\piprop_{i,j}(Y^j|X^i)
      \;
   \\[-1em]
   &
   \qquad\qquad\qquad\qquad\qquad
      \underbrace{
        \int_{Y^{\neg j}}
          \Big\{\prod_{k\neq j}^N \piprop_{i,k}(Y^k|X^i)\,\rmd Y^k\Big\} 
      }_{=1}
      \;\rmd Y^j
\\
   &
   =
   \frac{1}{N}\sum_{j=1}^N
      \alpha^{i,j}(X^i,z)\,\piprop_{i,j}(z|X^i) \;\rmd z
\end{align*}
because $\int_{Y^j} \delta_{Y^j}(\rmd z)\,dY^j=\rmd z$. The second term $A_2$ reads:
\begin{align*}
   A_2
   &
   =
   \int_{Y^{1:N}}
     \tilde\rho^i(X^i,Y)\,\delta_{X^i}(\rmd z)
     \;
     \Big\{\prod_{k=1}^N \piprop_{i,k}(Y^k|X^i)\,\rmd Y^k\Big\} 
\\
   &
   =
   \delta_{X^i}(\rmd z)\;
   \int_{Y^{1:N}}
     \Big\{1 - \frac{1}{N} \sum_{j=1}^N \alpha^{i,j}(X^i,Y^j)\Big\}
     \;
     \Big\{\prod_{k=1}^N \piprop_{i,k}(Y^k|X^i)\,\rmd Y^k\Big\} 
\\
   &
   =
   \delta_{X^i}(\rmd z)\;
   \Big\{
     1
     -
     \frac{1}{N} \sum_{j=1}^N 
       \int_{Y^{1:N}}
         \alpha^{i,j}(X^i,Y^j)
         \;
         \prod_{k=1}^N \piprop_{i,k}(Y^k|X^i)\,\rmd Y^k 
   \Big\}
\\
   &
   =
   \delta_{X^i}(\rmd z)\;
   \Big\{
     1
     -
     \frac{1}{N} \sum_{j=1}^N 
       \int_{Y^j}
         \alpha^{i,j}(X^i,Y^j)
         \;
         \piprop_{i,j}(Y^j|X^i) 
       \;\rmd Y^j
   \Big\}\,.
\end{align*}
Summing up $A_1$ and  $A_2$ proves the Lemma.
\carre

\subsection{Invariance property}

\begin{lemma}
\label{lemma.bilan.detaille.cond.MH1.alpha}
For all $(x,z)\in \R^n\times\R^n$ a.e. we have:
\begin{align*}
  \alpha^{i,j}(x,z) \;
  \pi(x) \;
  \piprop_{i,j}(z|x) 
  =
  \alpha^{i,j}(z,x)\;
  \pi(z)\;
  \piprop_{i,j}(x|z) 
  \,.
\end{align*}
\end{lemma}

\proof
For $(x,z)\not\in R^{i,j}$ the result is obvious. For $(x,z)\in R^{i,j}$ we have:
\begin{align*}
  &
  (r^{i,j}(x,z) \wedge 1) \;
  \pi(x) \;
  \piprop_{i,j}(z|x) 
  \\
  &
  \qquad\qquad\qquad
  =
  \min\!\Big\{
    \pi(z)\,
    \piprop_{i,j}(x|z) 
  \;,\;
    \pi(x)\,
    \piprop_{i,j}(z|x) 
  \Big\}
\\
  &
  \qquad\qquad\qquad
  =
  ( r^{i,j}(z,x) \wedge 1) \;
  \pi(z)\;
  \piprop_{i,j}(x|z) 
  \,.
\end{align*}
\carre

\bigskip

%
%

\begin{lemma}[conditional detailed balance]
\label{lemma.bilan.detaille.cond.MH1}
The following equality of measures defined on $\R^n \times\R^n$
\begin{align}
\label{eq.bilan.detaille.cond.MH1}
  P^i(Z^{1:i-1},X^{i:N};\rmd Z^i)
  \;
  \pi(X^i)\,\rmd X^i
  =
  P^i(Z^{1:i},X^{i+1:N};\rmd X^i)
  \;
  \pi(Z^i)\,\rmd Z^i
\end{align}
holds true for any  $i=1,\dots, N$, $Z^{1:i-1} \in \R^{(i-1)\times N}$, and $X^{i+1:N}\in\R^{(N-i) \times N}$.
\end{lemma}

\proof
Left hand side of  (\ref{eq.bilan.detaille.cond.MH1}) is a measure, say $\nu(\rmd Z^i\times \rmd X^i)$ on  $(\R^n\times\R^n,\BB(\R^n\times\R^n))$. For all  $A_1,\, A_2\in\BB(\R^n)$, we want to prove that $\nu(A_1\times A_2)=\nu(A_2\times A_1)$. We have: 
\begin{align*}
  \nu(A_1\times A_2)
  &=
  \int   
    P^i(Z^{1:i-1},X^{i:N};A_1) \,
    \indic_{A_2}(X^i)\,  
    \pi(X^i)\,
    \rmd X^i
\end{align*}
and
\begin{align*}
  P^i(Z^{1:i-1},X^{i:N};A_1)
  &=
  \frac{1}{N}\sum_{j=1}^N 
  \int 
      \indic_{A_1}(Z^i)\,
      \alpha^{i,j}(X^i,Z^i)
      \;
      \piprop_{i,j}(Z^i|X^i)
      \;
      \rmd Z^i
  \\[-0.7em]
  &\qquad\qquad\qquad\qquad
   +
   \rho^i(X^i) \;
   \indic_{A_1}(X^i)
\end{align*}
so that 
\begin{align}
\nonumber
  &
  \nu(A_1\times A_2)
  \\
\nonumber
  &
  \qquad
  =
  \frac{1}{N}\sum_{j=1}^N 
  \iint
      \indic_{A_1}( Z^i) \,
      \indic_{A_2}( X^i)\,  
      \alpha^{i,j}( X^i, Z^i) \,
      \pi( X^i)\,
      \piprop_{i,j}( Z^i|X^i)  \,
      \rmd  X^i\,
      \rmd Z^i
  \\
\label{eq.lemma.bilan.detaille.cond.MH1}
  &
  \qquad\qquad
  +
  \int   
    \rho^i( X^i) \,
    \indic_{A_1}( X^i) \,
    \indic_{A_2}( X^i)\,  
    \pi( X^i)\,
    \rmd  X^i
    \,.
\end{align}
And from Lemma \ref{lemma.bilan.detaille.cond.MH1.alpha}, we get:
\begin{align*}
  &
  \nu(A_1\times A_2)
\\
  &
  \qquad
  =
  \frac{1}{N}\sum_{j=1}^N 
  \iint
      \indic_{A_1}( Z^i)\,
      \indic_{A_2}( X^i)\;
      \alpha^{i,j}( Z^i, X^i)  \;
      \pi( Z^i)\,
      \piprop_{i,j}( X^i|Z^i) \,
      \rmd  Z^i\,
      \rmd  X^i
  \\
  &\qquad\qquad
  +
  \int   
    \rho^i( X^i) \,
    \indic_{A_1}(  X^i) \,
    \indic_{A_2}(  X^i)\,  
    \pi( X^i)\,
    \rmd   X^i
\end{align*}
Exchanging the name of variables $ X^i\leftrightarrow  Z^i$ in the first term of the right hand side of the previous equality, leads to the same expression as  (\ref{eq.lemma.bilan.detaille.cond.MH1}) where $A_1$ and $A_2$ were interchanged, in other words   $\nu(A_1\times A_2) =  \nu(A_2\times A_1)$.
\carre

\bigskip

\begin{proposition}[invariance]
The probability measure
\[
 \Pi({\rmd} X)=\pi(X^1)\,\rmd X^1 \cdots \pi(X^N)\,\rmd X^N
\]
is an invariant distribution of the Markov kernel  $P$, i.e.  $\Pi P = \Pi$ that is:
\begin{align}
\label{eq.invariance.MHG.MH1}
 \int_X  P(X,\rmd Z) \; \Big\{\prod_{i=1}^N \pi(X^i)\,\rmd X^i\Big\}
 =
 \prod_{i=1}^N \pi(Z^i)\,\rmd Z^i
 \,.
\end{align}
\end{proposition}

\proof
\begin{align*}
 &
 \int_X  
   P(X,\rmd Z) \; 
   \Big\{ \prod_{i=1}^N \pi(X^i)\,\rmd X^i \Big\}
\\
 &\qquad\qquad
 = 
 \int_X  
    P^1(X^{1:N};\rmd Z^{1}) \; 
    P^2(Z^{1},X^{2:N};\rmd Z^{2}) \cdots 
  \\[-0.8em]
  &\qquad\qquad\qquad\qquad\qquad\qquad
    \cdots
    P^N(Z^{1:N-1},X^{N};\rmd Z^{N}) \;
    \Big\{ \prod_{i=1}^N \pi(X^i)\,\rmd X^i \Big\}
\\
 &\qquad\qquad
 = 
 \int_X
   P^1(X^{1:N};\rmd Z^{1})  \;
     \pi(X^1) \, \rmd X^1\;
   P^2(Z^{1},X^{2:N};\rmd Z^{2}) \cdots 
  \\[-0.8em]
  &\qquad\qquad\qquad\qquad\qquad\qquad
    \cdots
   P^n(Z^{1:N-1},X^{N};\rmd Z^{N})\;
   \Big\{\prod_{i=2}^N \pi(X^i) \, \rmd X^i\Big\}
   \,.
\end{align*}
Using (\ref{eq.bilan.detaille.cond.MH1}) with  $i=1$ gives:
\begin{align*}
 &
 \int_X  
   P(X,\rmd Z) \; 
   \Big\{ \prod_{i=1}^N \pi(X^i)\,\rmd X^i \Big\}
 =
\\
 &\qquad
 = 
 \int_X 
   P^1(Z^1,X^{2:N};\rmd X^{1})\;
    \pi(Z^1) \, \rmd Z^1
   \; 
   P^2(Z^{1},X^{2:N};\rmd Z^{2}) \cdots 
  \\[-0.8em]
  &\qquad\qquad\qquad\qquad\qquad\qquad
    \cdots
   P^n(Z^{1:N-1},X^N;\rmd Z^N) \;
   \Big\{ 
     \prod_{i=2}^N \pi(X^i) \, \rmd X^i 
   \Big\}
   \,.
\end{align*}
In this last expression the kernel $P^1(Z^1,X^{2:N};\rmd X^{1})$ is a measure on the variable $X^{1}$ which no longer appears in the integrand. Therefore its integral with respect to this variable is  1, hence:
\begin{align*}
 &
 \int_X  P(X,\rmd Z) \; \Big\{ \prod_{i=1}^N \pi(X^i)\,\rmd X^i\Big\}
 =
\\
 &\qquad 
 = 
 \pi(Z^1) \, \rmd Z^1\;
 \int_{X^{2:N}}
    P^2(Z^{1},X^{2:N};\rmd Z^{2}) \cdots 
  \\[-0.8em]
  &\qquad\qquad\qquad\qquad\qquad\qquad
    \cdots
    P^n(Z^{1:N-1},X^N;\rmd Z^N) \;
    \Big\{
      \prod_{i=2}^N \pi(X^i) \, \rmd X^i
    \Big\}
    \,.
\end{align*}
Repeating this procedure successively for  $X^2$ to $X^N$ leads to (\ref{eq.invariance.MHG.MH1}).
\carre

\section[Parallel/interacting MwG algorithm]
        {Parallel/interacting Metropolis within Gibbs\\ (MwG) algorithm}
\label{sec.para.int.MwG}

Let $\pi(x)$ be the probability density function of a target distribution defined on $(\R^n,\BB(\R^n))$. For  $\ell=1,\dots, n$, we define the conditional laws:
\begin{align}
\label{eq.loi.cond}
  \pi_\ell(x_\ell|x_{\neg \ell})
  \eqdef
  \frac{\pi(x_{1:n})}{\int \pi(x_{1:n}) \,\rmd x_{\neg\ell}}
  \,.
\end{align}
When we know to sample from  (\ref{eq.loi.cond}), we are able to use the Gibbs sampler. It is possible to adapt our interacting method to parallel Gibbs sampler. But very often we do not know how to sample from (\ref{eq.loi.cond}) and therefore we consider proposal conditional densities $\piprop_\ell(x_\ell)$ defined for all $\ell$. In this case, we use Metropolis within Gibbs algorithm (see appendix). We present in the following how to make interactions between parallel MwG algorthims. The MwG algorithm is more general than Gibbs algorithm, so a parallel/interacted Gibbs algorithm can easily be deduced from the parallel/interacted MwG algorithm.  

\subsection{The  algorithm}
\label{sec.MwG.algo}

One iteration $X\to Z$ of the parallel/interacting Metropolis within Gibbs method consists in updating  the components $X_\ell$ successively for $\ell=1,\dots,n$, i.e.
\[
  [X_{1:n}]
  \to
  [Z_1 X_{2:n}]
  \to
  [Z_{1:2} X_{3:n}]
  \cdots
  [Z_{1:n-1} X_{n}]
  \to
  [Z_{1:n}]
  \,.
\]
For each $\ell$ fixed, the subcomponents  $X^i_\ell$ are updated sequentially for $i=1,\dots,N$ in two steps:
\begin{enumerate}

\item
{\it  Proposal step:}
We sample independently $N$ candidates  $Y_\ell^{j}\in\R$ for  $j=1:N$ according to:
\[
     Y^j_\ell 
     \sim \pipropl_{i,j}(\xi|\bbra{Z,X^i_\ell,X}^i_\ell)\,\rmd \xi
     \,,
     \qquad 1\leq j\leq n
\]
where
\[
   \bbra{Z,\xi,X}^i_\ell
   \eqdef 
   \left[
     Z_{1:\ell-1} 
     \left|
     \begin{smallmatrix}
         Z^1_{\ell} \\
         \vdots\\
         Z^{i-1}_{\ell} \\
         \xi \\
         X^{i+1}_{\ell} \\
         \vdots\\
         X^{N}_{\ell} \\
     \end{smallmatrix}
     \right|
     X_{\ell+1:n} 
   \right]
   \,.
\]
We also use the following lighter notation:
\[
   \pipropl_{i,j}(\xi|\xi')
   =
   \pipropl_{i,j}(\xi|\bbra{Z,\xi',X}^i_\ell)
   \,.
\]


\item
{\it Selection step:}
The subcomponent $X^i_\ell$ could be replaced by one of the $N$ candidates $Y^{1:N}_\ell$ or stay unchanged according to a multinomial sampling, the resulting value is called $Z^i_\ell$, i.e.:
\[
  Z^i_\ell
  \ot
  \left\{\begin{array}{ll}
     Y^1_\ell 
     & \textrm{with probability }
       \frac{1}{N}\,\alpha^{i,1}_\ell(X^i_\ell,Y^1_\ell)
       \,,
     \\
     \ \vdots
     \\
     Y^N_\ell 
     & \textrm{with probability }
       \frac{1}{N}\,\alpha^{i,N}_\ell(X^i_\ell,Y^N_\ell)
       \,,
     \\[0.5em]
     X^i_\ell 
     & \textrm{with probability }
       \tilde\rho^i_\ell(X^i_\ell,Y^{1:N}_\ell)
  \end{array}\right.
\]
where:
\begin{align*}
  \alpha^{i,j}_\ell(\xi,\xi')
  &\eqdef
  \frac{\pi_\ell(\xi'|X^i_{\neg\ell})}
       {\pi_\ell(\xi |X^i_{\neg\ell})}
  \,
  \frac{\pipropl_{i,j}(\xi |\xi')}
       {\pipropl_{i,j}(\xi'|\xi )}
  \wedge 1
  \,,
\\
  \tilde\rho^i_\ell(X^i_\ell,Y^{1:N}_\ell)
  &\eqdef
  1-\frac{1}{N}\sum_{j=1}^N\alpha^{i,j}_\ell(X^i_\ell,Y^j_\ell)
  \,.
\end{align*}

\end{enumerate}
The resulting algorithm is depicted in Algorithm~\ref{algo.mh.gibbs.1}.

\begin{algorithm}[ht]
\caption{\sl Parallel/interacting MwG.}
\label{algo.mh.gibbs.1}
\begin{center}
\begin{minipage}{12cm}
\hrulefill\\[-1em]
\mbox{}
\begin{algorithmic}
\STATE choose $X\in\R^{n\times N}$  
\FOR {$k=1,2,\dots$}
  \FOR {$\ell=1:n$}
    \FOR {$i=1:N$}
      \FOR {$j=1:N$}
         \STATE $Y^j_\ell\sim \pipropl_{i,j}(\xi)\,\rmd \xi$  
         \STATE $\alpha^j
            \ot 
            \frac{\pi_\ell(Y^j_\ell|X^i_{\neg\ell})}
                 {\pi_\ell(X^i_\ell|X^i_{\neg\ell})}
             \,
            \frac{\pipropl_{i,j}(X^i_\ell|Y^j_\ell)}
                 {\pipropl_{i,j}(Y^j_\ell|X^i_\ell)}
                  \wedge 1
                      $
      \ENDFOR
      \STATE $\tilde\rho \ot 1 - \frac{1}{N}\sum_{j=1}^N\alpha^j$
      \STATE $X^i_\ell
        \ot
        \left\{\begin{array}{ll}
          Y^1_\ell &\textrm{with probability } \alpha^1/N
          \\[-0.4em]
          \ \vdots
          \\
          Y^{N}_\ell &\textrm{with probability } \alpha^N/N
          \\[0.2em]
          X^i_\ell
          &\textrm{with probability } 
             \tilde\rho
     \end{array}\right.$
    \ENDFOR
  \ENDFOR
  \ENDFOR
\end{algorithmic}
\hrulefill
\end{minipage}
\end{center}
\end{algorithm}

\subsection{Description of the MH kernel}

\begin{lemma}
The Markov kernel on $\R^{n\times N}$  associated with the MH algorithm described in Section \ref{sec.MwG.algo}, is
\begin{align}
\label{eq.P}
  P(X,\rmd Z)
  \eqdef
  P_1(X_{1:n};\rmd Z_{1})
  \;
  P_2(Z_{1},X_{2:n};\rmd Z_{2})
  \cdots
  P_n(Z_{1:n-1},X_{n};\rmd Z_{n})
  \,.
\end{align}
At iteration $\ell$, the kernel $P_\ell(Z_{1:\ell-1},X_{\ell:n};\rmd Z_{\ell})$ generates $Z^{1:N}_{\ell}$ from the already updated components $Z^{1:N}_{1:\ell-1}$ and the remaining components  $X^{1:N}_{\ell:n}$. 

Each component  $Z^i_{1:\ell}$, for  $i=1\cdots N$, is updated independently one from each other:
\begin{align}
\label{eq.P.ell}
  &
  P_\ell(Z_{1:\ell-1},X_{\ell:n};\rmd Z_\ell)  
  \eqdef
  \prod_{i=1}^N
    P_\ell^i(\bbra{Z,X_\ell^i,X}^i_\ell;\rmd Z_\ell^i) 
  \,. 
\end{align}
Here $Z_\ell^i$ is generated from $\bbra{Z,X_\ell^i,X}^i_\ell$ according to:
\begin{align}
\label{eq.P.ell.i.sol}
  P_\ell^i(\bbra{Z,\xi,X}^i_\ell;\rmd\xi')  
  \eqdef
  \frac{1}{N}\sum_{j=1}^N 
      \alpha^{i,j}_\ell(\xi,\xi')
      \;
      \pipropl_{i,j}(\xi'|\xi)
      \;
      \rmd\xi'
  +
  \rho^i_\ell(\xi)
  \;
  \delta_\xi(\rmd\xi')
\end{align}
Acceptation probabilities are:
\begin{align}
\label{eq.alpha}
  \alpha^{i,j}_\ell(\xi,\xi')
  &\eqdef
  \left\{\begin{array}{ll}
     r^{i,j}_\ell(\xi,\xi') \wedge 1
     & 
     \textrm{ if }(\xi,\xi') \in R^{i,j}_\ell\,,
   \\[0.3em]
     0
     & 
     \textrm{ otherwise,}
  \end{array}\right.
\\[0.8em]
\label{eq.r}
  r^{i,j}_\ell(\xi,\xi')
  &\eqdef
  \frac{\pi_\ell(\xi'|Z_{1:\ell-1}^i,X_{\ell+1:n}^i)}
       {\pi_\ell(\xi |Z_{1:\ell-1}^i,X_{\ell+1:n}^i)}
  \,
  \frac{\pipropl_{i,j}(\xi |\xi')}
       {\pipropl_{i,j}(\xi'|\xi )}
   \,,
\\[0.8em]
\label{eq.r.ell.i}
  \rho^i_\ell(\xi)
  &\eqdef
  1
  -
  \frac{1}{N}
   \sum_{j=1}^N  \int_\R 
     \alpha^{i,j}_\ell(\xi,\xi') \;
     \pipropl_{i,j}(\xi'|\xi) \;
     \rmd \xi'
   \,.
\end{align}
Finally, $R^{i,j}_\ell$ is the set of ordered pairs $(\xi,\xi')\in\R^2$ such that
\[
  \begin{array}{r}
    \pi_\ell(\xi'|Z_{1:\ell-1}^i,X_{\ell+1:n}^i)  
       \;
       \pipropl_{i,j}(\xi |\xi') > 0  \,,
   \\[0.3em]
    \pi_\ell(\xi |Z_{1:\ell-1}^i,X_{\ell+1:n}^i)    
       \;
       \pipropl_{i,j}(\xi'|\xi )  >0    \,.
  \end{array}
\]
Note that the functions  $\alpha^{i,j}_\ell(\xi,\xi')$, $\rho^i_\ell(\xi)$, $r^{i,j}_\ell(\xi,\xi')$ and the set  $R^{i,j}_\ell$ depend on  $Z_{1:\ell-1}$ and $X_{\ell+1:n}$.
\end{lemma}

\proof
This construction follows the general setup proposed by Luke Tierney in \cite{tierney1998a}.
The kernel is defined by: 
\begin{align*}
  P_\ell^i(\bbra{Z,\xi,X}^i_\ell;\rmd \xi')  
  &\eqdef
  \int_{\R^N}
    \underbrace{
      S^i_\ell(\bbra{Z,\xi,X}^i_\ell,\zeta^{1:N};\rmd \xi')
    }_{\textrm{\tiny selection kernel}}
    \times
    \underbrace{
      Q^i_\ell(\bbra{Z,\xi,X}^i_\ell;\rmd \zeta^{1:N})
    }_{\textrm{\tiny proposal kernel}}
    \,.
\end{align*}
This kernel consists firstly in proposing a population of $N$ candidates $\zeta^{1:N}\in\R^N$ sampled from:
\begin{align}
\label{eq.Q.ell}
  Q^i_\ell(\bbra{Z,\xi,X}^i_\ell;\rmd \zeta^{1:N})
  &
  \eqdef
  \prod_{j=1}^N 
     \pipropl_{i,j}(\zeta ^j|\xi)\; 
     \rmd \zeta ^j
 \,,
\end{align}
then secondly in selecting among these candidates or rejecting them according to a MH technique, i.e.
\begin{align}
\label{eq.S.ell.i}
  S^i_\ell(\bbra{Z,\xi,X}^i_\ell,\zeta^{1:N};\rmd \xi')
  &
  \eqdef
  \frac{1}{N} \sum_{j=1}^N
    \alpha^{i,j}_\ell(\xi,\zeta^j) \, 
    \delta_{\zeta^j}(\rmd \xi')
  +
  \tilde\rho^i_\ell(\xi,\zeta^{1:N}) \, \delta_{\xi}(\rmd \xi')
\end{align}
where $\alpha^{i,j}_\ell$ is given by  (\ref{eq.alpha}) and $\tilde\rho^i_\ell(\xi,\zeta^{1:N})
  \eqdef
  1 -  \frac{1}{N} \sum_{j=1}^N \alpha^{i,j}_\ell(\xi,\zeta^j)$.

Hence:  
\begin{align*}
  &
  P_\ell^i(\bbra{Z,\xi,X}^i_\ell;\rmd \xi')  
  \eqdef
  \int_{\zeta^{1:N}}
      S^i_\ell(\bbra{Z,\xi,X}^i_\ell,\zeta^{1:N};\rmd \xi')
    \;
      Q^i_\ell(\bbra{Z,\xi,X}^i_\ell;\rmd \zeta^{1:N})
\\ 
  &
  \qquad
  =
  \frac{1}{N} \sum_{j=1}^N 
  \int_{\zeta^{1:N}}
    \alpha^{i,j}_\ell(\xi,\zeta^j) \, 
    \delta_{\zeta^j}(\rmd \xi')\,
    \prod_{k=1}^N 
     \pipropl_{i,j}(\zeta ^k|\xi)\; 
     \rmd \zeta ^k
  \\
  &
  \qquad\qquad\qquad
  +
  \Big\{
  1
  -
  \frac{1}{N} \sum_{j=1}^N 
  \int_{\zeta^{1:N}}
    \alpha^{i,j}_\ell(\xi,\zeta^j) \, 
    \prod_{k=1}^N 
     \pipropl_{i,j}(\zeta ^k|\xi)\; 
     \rmd \zeta ^k
  \Big\} \; \delta_{\xi}(\rmd \xi')
\\ 
  &
  \qquad
  =
  \frac{1}{N} \sum_{j=1}^N 
  \int_{\zeta^j}
    \alpha^{i,j}_\ell(\xi,\zeta^j) \, 
    \delta_{\zeta^j}(\rmd \xi')\,
     \piprop_\ell(\zeta ^j|\xi)\; 
     \rmd \zeta ^j
  \\
  &
  \qquad\qquad\qquad
  +
  \Big\{
  1
  -
  \frac{1}{N} \sum_{j=1}^N 
  \int_{\zeta^j}
    \alpha^{i,j}_\ell(\xi,\zeta^j) \, 
     \piprop_\ell(\zeta ^j|\xi)\; 
     \rmd \zeta ^j
  \Big\} \; \delta_{\xi}(\rmd \xi')
\\ 
  &
  \qquad
  =
  \frac{1}{N} \sum_{j=1}^N 
    \alpha^{i,j}_\ell(\xi,\xi') \, 
     \piprop_\ell(\xi'|\xi)\; 
     \rmd \xi'
  \\
  &
  \qquad\qquad\qquad
  +
  \Big\{
  1
  -
  \frac{1}{N} \sum_{j=1}^N 
  \int_{\xi''}
    \alpha^{i,j}_\ell(\xi,\xi'') \, 
     \piprop_\ell(\xi''|\xi)\; 
     \rmd \xi''
  \Big\} \; \delta_{\xi}(\rmd \xi')
\end{align*}
which correspond to Equations (\ref{eq.P.ell.i.sol}) to (\ref{eq.r.ell.i}).
\carre

\subsection{Invariance property}

\begin{lemma}
\label{lemma.alpha}
For almost all  $(\xi,\xi')\in\R^2$:
\begin{align*}
  &
  \alpha^{i,j}_\ell(\xi,\xi') \;
  \pi_\ell(\xi|Z_{1:\ell-1}^i,X_{\ell+1:n}^i) \;
  \pipropl_{i,j}(\xi'|\xi) 
 \\
  &\qquad\qquad\qquad\qquad\qquad
  =
  \alpha^{i,j}_\ell(\xi',\xi)  \;
  \pi_\ell(\xi'|Z_{1:\ell-1}^i,X_{\ell+1:n}^i)\;
  \pipropl_{i,j}( \xi|\xi') 
\end{align*}
for any $\ell$, $i$, $j$, $(Z_{1:\ell-1}^i,X_{\ell+1:n}^i)$, and
$(Z_{1:\ell-1}^j,X_{\ell+1:n}^j)$.
\end{lemma}

\proof
For $(\xi,\xi')\not\in R^{i,j}_\ell$, the result is obvious.
For $(\xi,\xi')\in R^{i,j}_\ell$ a.e.:
\begin{align*}
  &
  (r^{i,j}_\ell(\xi,\xi') \wedge 1) \;
  \pi_\ell(\xi|Z_{1:\ell-1}^i,X_{\ell+1:n}^i) \;
  \piprop_\ell(\xi'|\xi) 
 \\
  &\qquad
  =
  \min\!\Big\{
    \pi_\ell(\xi'|Z_{1:\ell-1}^i,X_{\ell+1:n}^i)\,
    \piprop_\ell( \xi|\xi') 
  \;,\;
    \pi_\ell(\xi|Z_{1:\ell-1}^i,X_{\ell+1:n}^i)\,
    \piprop_\ell(\xi'|\xi) 
  \Big\}
\\
  &\qquad
  =
  ( r^{i,j}_\ell(\xi',\xi) \wedge 1) \;
  \pi_\ell(\xi'|Z_{1:\ell-1}^i,X_{\ell+1:n}^i)\;
  \piprop_\ell( \xi|\xi') 
  \,.
\end{align*}
\carre

\bigskip

\begin{lemma}[conditional detailed balance]
\label{lemma.bilan.detaille.cond}
The following e\-qua\-lity of measures defined on $\R \times\R$
\begin{align}
\nonumber
  &
  P_\ell^i(\bbra{Z,\xi,X}^i_\ell;\rmd\xi')
  \times
  \pi_\ell(\xi|Z_{1:\ell-1}^i,X_{\ell+1:n}^i)\,\rmd \xi
  \\
\label{eq.bilan.detaille.cond}
  & \qquad\qquad\qquad
  =
  P_\ell^i(\bbra{Z,\xi',X}^i_\ell;\rmd \xi)
  \times
  \pi_\ell(\xi'|Z_{1:\ell-1}^i,X_{\ell+1:n}^i)\,\rmd\xi'
\end{align}
holds true for any  $\ell=1\cdots n$, $i=1\cdots N$ and $Z_{1:\ell-1} \in \R^{N\times(\ell-1)}$, $X_{\ell+1:n}\in\R^{N\times (n-\ell)}$.
\end{lemma}

\proof
The left hand side of equality (\ref{eq.bilan.detaille.cond}) is a measure  $\nu(\rmd\xi'\times\rmd \xi)$ defined on  $(\R^2,\BB(\R^2))$. For all $A_1,\, A_2\in\BB(\R)$, we want to prove that $\nu(A_1\times A_2)=\nu(A_2\times A_1)$.

We have: 
\begin{align*}
  \nu(A_1\times A_2)
  &=
  \int   
    P_\ell^i(\bbra{Z,\xi,X}^i_\ell;A_1) \,
    \indic_{A_2}(\xi)\,  
    \pi_\ell(\xi|Z_{1:\ell-1}^i,X_{\ell+1:n}^i)\,
    \rmd \xi
\end{align*}
and
\begin{align*}
  &
  P_\ell^i(\bbra{Z,\xi,X}^i_\ell;A_1)
  =
  \frac{1}{N}\sum_{j=1}^N 
  \int 
      \indic_{A_1}(\xi')\,
      \alpha^{i,j}_\ell(\xi,\xi')
      \;
      \piprop_\ell(\xi'|\xi)
      \;
      \rmd\xi'
   +
   \rho^i_\ell(\xi) \;
   \indic_{A_1}(\xi)
\end{align*}
so that
\begin{align}
\nonumber
  \nu(A_1\times A_2)
  &=
  \frac{1}{N}\sum_{j=1}^N 
  \iint
      \indic_{A_1}(\xi') \,
      \indic_{A_2}(\xi)\,  
      \alpha^{i,j}_\ell(\xi,\xi') \,
  \\[-0.8em]
\nonumber
  &\qquad\qquad\qquad\qquad      
      \pi_\ell(\xi|Z_{1:\ell-1}^i,X_{\ell+1:n}^i)\,
      \piprop_\ell(\xi'|\xi)  \,
      \rmd \xi\,
      \rmd\xi'
  \\
\label{eq.lemma.bilan.detaille.cond}
  &\qquad\qquad
  +
  \int   
    \rho^i_\ell(\xi) \,
    \indic_{A_1}(\xi) \,
    \indic_{A_2}(\xi)\,  
    \pi_\ell(\xi|Z_{1:\ell-1}^i,X_{\ell+1:n}^i)\,
    \rmd \xi
\end{align}
Using Lemma \ref{lemma.alpha} we get:
\begin{align*}
  \nu(A_1\times A_2)
  &=
  \frac{1}{N}\sum_{j=1}^N 
  \iint
      \indic_{A_1}(\xi')\,
      \indic_{A_2}(\xi)\;
      \alpha^{i,j}_\ell(\xi',\xi) 
  \\[-0.8em]
\nonumber
  &\qquad\qquad\qquad \qquad      
      \pi_\ell(\xi'|Z_{1:\ell-1}^i,X_{\ell+1:n}^i)\,
      \piprop_\ell(\xi|\xi') \,
      \rmd \xi'\,
      \rmd \xi
  \\
  &\qquad\qquad
  +
  \int   
    \rho^i_\ell(\xi) \,
    \indic_{A_1}( \xi) \,
    \indic_{A_2}( \xi)\,  
    \pi_\ell(\xi|Z_{1:\ell-1}^i,X_{\ell+1:n}^i)\,
    \rmd  \xi
\end{align*}
Exchanging the name of variables $\xi\leftrightarrow \xi'$ in the first term of the right hand side of the previous equality leads to the same expression as  (\ref{eq.lemma.bilan.detaille.cond})  where $A_1$ and $A_2$ were interchanged, in other words   $\nu(A_1\times A_2) =  \nu(A_2\times A_1)$.
\carre

\bigskip

\begin{proposition}[invariance]
The measure
\[
 \Pi({\rmd} X)=\pi(X^1)\,\rmd X^1 \cdots \pi(X^N)\,\rmd X^N
\]
is invariant for the kernel $P$, that is $\Pi P = \Pi$ i.e.:
\begin{align}
\label{eq.invariance.MHG}
 \int_X  P(X,\rmd Z) \; \Big\{\prod_{i=1}^N \pi(X^i)\,\rmd X^i\Big\}
 =
 \prod_{i=1}^N \pi(Z^i)\,\rmd Z^i
 \,.
\end{align}
\end{proposition}

\proof
\begin{align*}
 &
 \int_X  
   P(X,\rmd Z) \; 
   \Big\{ \prod_{i=1}^N \pi(X^i)\,\rmd X^i \Big\}
\\
 &\quad
 = 
 \int_X  
    P_1(X_{1:n};\rmd Z_{1}) \; 
    P_2(Z_{1},X_{2:n};\rmd Z_{2}) \cdots 
    P_n(Z_{1:n-1},X_{n};\rmd Z_{n}) \;
  \\[-0.8em]
  &\qquad\qquad\quad
    \prod_{i=1}^N 
    \big\{
       \pi_1(X_1^i|X_{2:n}^i)  \, \rmd X_1^i
       \;
       \pi_{\neg 1}(X^i_{2:n}) \, \rmd X^i_{2:n}
    \big\}
\\
 &\quad
 = 
 \int_X
   P_1(X_{1:n};\rmd Z_{1})  \;
     \Big\{\prod_{i=1}^N \pi_1(X_1^i|X_{2:n}^i) \, \rmd X_1^i\Big\}\;
  \\[-0.8em]
  &\qquad\qquad\quad
   P_2(Z_{1},X_{2:n};\rmd Z_{2})    \cdots
   P_n(Z_{1:n-1},X_{n};\rmd Z_{n})\;
   \Big\{\prod_{i=1}^N \pi_{\neg 1}(X^i_{2:n}) \, \rmd X^i_{2:n}\Big\}
\\
 &\quad
 = 
 \int_X
     \Big\{\prod_{i=1}^N P^i_1(\bbra{Z,X^i_1,X}_1^i;dZ^i_1) \Big\}
     \;
     \Big\{\prod_{i=1}^N \pi_1(X_1^i|X_{2:n}^i) \, \rmd X_1^i\Big\}\;
  \\[-0.8em]
  &\qquad\qquad\quad
   P_2(Z_{1},X_{2:n};\rmd Z_{2}) \cdots 
   P_n(Z_{1:n-1},X_{n};\rmd Z_{n})\;
   \Big\{\prod_{i=1}^N \pi_{\neg 1}(X^i_{2:n}) \, \rmd X^i_{2:n}\Big\}
\end{align*}
Moreover
\begin{align*}
  &
  P_1(X_{1:n};\rmd Z_{1})  \;
  \Big\{\prod_{i=1}^N \pi_1(X_1^i|X_{2:n}^i) \, \rmd X_1^i\Big\}
  =
\\
  &\qquad
  =
  \Big\{\prod_{i=1}^N P^i_1(\bbra{Z,X^i_1,X}_1^i;dZ^i_1) \Big\}
  \;
  \Big\{\prod_{i=1}^N \pi_1(X_1^i|X_{2:n}^i) \, \rmd X_1^i\Big\}
\\
  &\qquad
  =
  \prod_{i=1}^N 
    P^i_1(\bbra{Z,X^i_1,X}_1^i;dZ^i_1) 
    \;
    \pi_1(X_1^i|X_{2:n}^i) \, \rmd X_1^i
\\
  &\qquad
  =
  \prod_{i=1}^N 
    P^i_1(\bbra{Z,Z^i_1,X}_1^i;dX^i_1) 
    \;
    \pi_1(Z_1^i|X_{2:n}^i) \, \rmd Z_1^i
\end{align*}
this last equality follows from Equation (\ref{eq.bilan.detaille.cond}). Hence,
\begin{align*}
 &
 \int_X  
   P(X,\rmd Z) \; 
   \Big\{ \prod_{i=1}^N \pi(X^i)\,\rmd X^i \Big\}
\\
 &\quad
 = 
 \int_X
  \prod_{i=1}^N 
     \Big\{
       P^i_1(\bbra{Z,Z^i_1,X}_1^i;dX^i_1) 
      \;
      \pi_1(Z_1^i|X_{2:n}^i) \, \rmd Z_1^i
     \Big\}
     \;
   P_2(Z_{1},X_{2:n};\rmd Z_{2}) \cdots 
  \\[-0.8em]
  &\qquad\qquad\qquad\qquad
    \cdots
   P_n(Z_{1:n-1},X_{n};\rmd Z_{n})\;
   \Big\{\prod_{i=1}^N \pi_{\neg 1}(X^i_{2:n}) \, \rmd X^i_{2:n}\Big\}
\end{align*}
In this last expression, for $i=1,\dots,N$, the kernel $P^i_1(\bbra{Z,Z^i_1,X}_1^i;dX^i_1)$ is a measure for the variable $X^i_{1}$ which no longer appears in the integrand. Using the fact that the integral of the kernel w.r.t. $X^i_{1}$ is 1 we get:

\begin{align*}
 &
 \int_X  
   P(X,\rmd Z) \; 
   \Big\{ \prod_{i=1}^N \pi(X^i)\,\rmd X^i \Big\}
\\
 &\quad
 = 
 \int_{X_{2:N}}
  \prod_{i=1}^N 
     \Big\{
      \pi_1(Z_1^i|X_{2:n}^i) \, \rmd Z_1^i
     \Big\}
     \;
   P_2(Z_{1},X_{2:n};\rmd Z_{2}) \cdots 
  \\[-0.8em]
  &\qquad\qquad\qquad\qquad
    \cdots
   P_n(Z_{1:n-1},X_{n};\rmd Z_{n})\;
   \Big\{\prod_{i=1}^N \pi_{\neg 1}(X^i_{2:n}) \, \rmd X^i_{2:n}\Big\}
\\
 &\quad
 = 
 \int_{X_{2:N}}
  \prod_{i=1}^N 
   P_2(Z_{1},X_{2:n};\rmd Z_{2}) \cdots 
  \\[-0.8em]
  &\qquad\qquad\qquad\qquad
    \cdots
   P_n(Z_{1:n-1},X_{n};\rmd Z_{n})\;
   \Big\{\prod_{i=1}^N \pi(Z_1^iX^i_{2:n}) \, \rmd Z_1^i\,\rmd X^i_{2:n}\Big\}
\end{align*}

Repeating this process successively for $X_2$ to $X_n$ leads to 
(\ref{eq.invariance.MHG}).
\carre

\section{Numerical tests}
\label{sec.tests}

\subsection{A multi-modal example}
\label{sec.multi-modal}

\begin{figure}[p]
\centering
\includegraphics[width=2.5in]{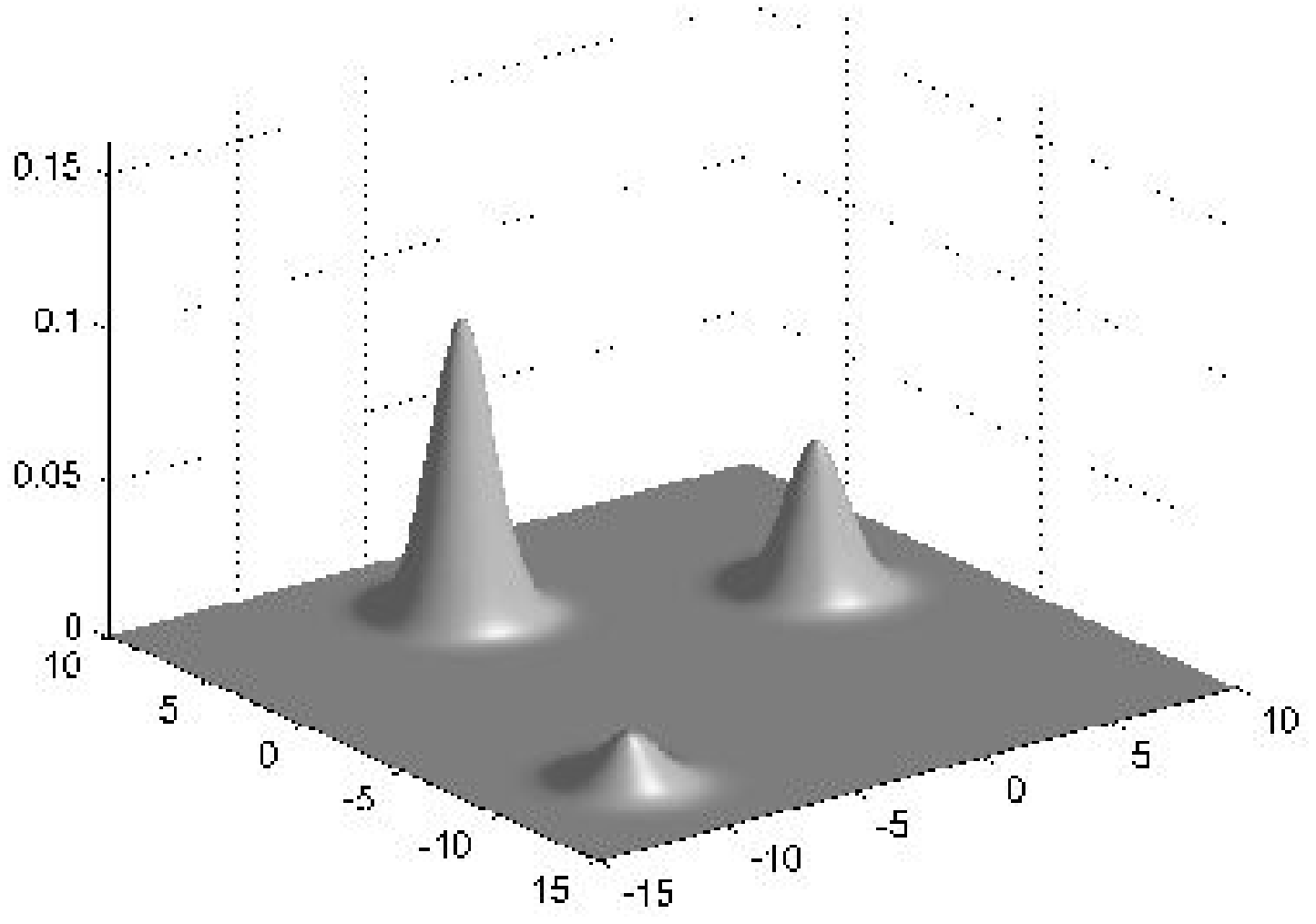}
\includegraphics[width=2.5in]{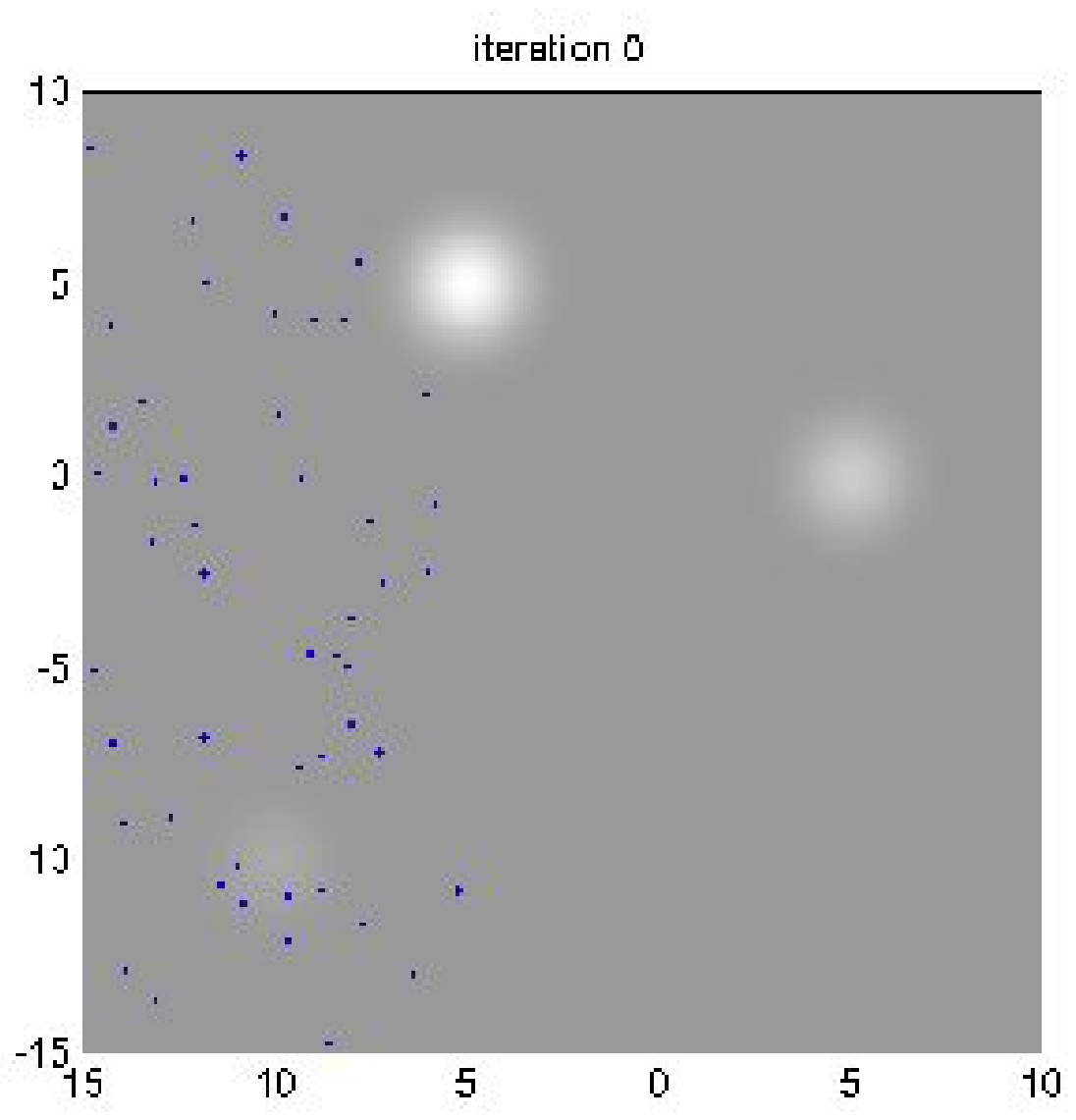}
\caption{Target distribution $\pi(x)$ (left) and  initial positions of the chains $\sfX^{(0),i}$, for $i=1\cdots N$ (right).}
\label{graphe.multi.target}
\vskip2em
\includegraphics[width=2.5in]{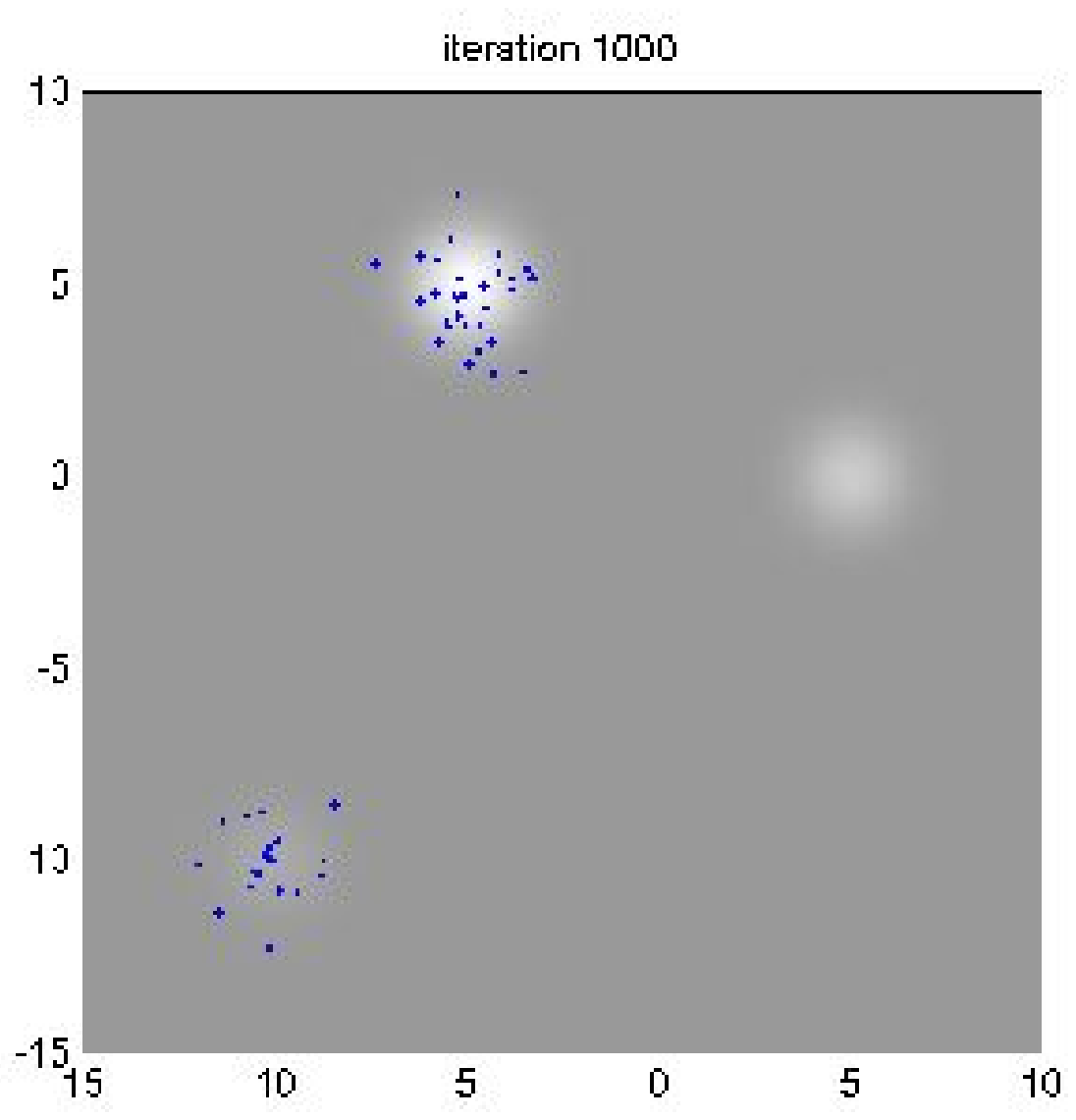}
\includegraphics[width=2.5in]{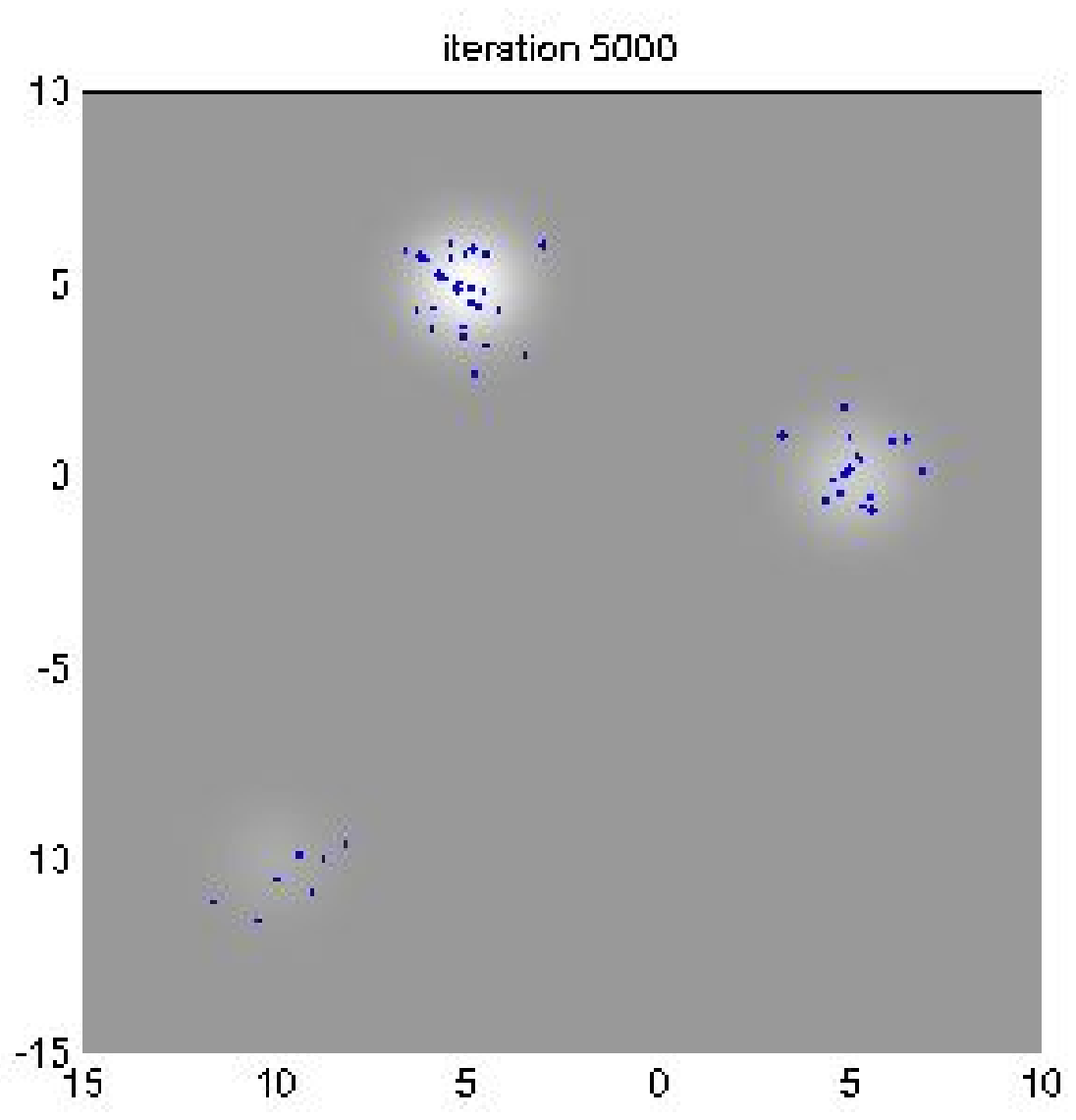}
\caption{Positions of the chains $\sfX^{(k),i}$, for $i=1\cdots N$, at iterations $k=1000$ (left) and $k=5000$  (right).}
\label{graphe.multi.target2}
\vskip2em
\includegraphics[width=2.5in]{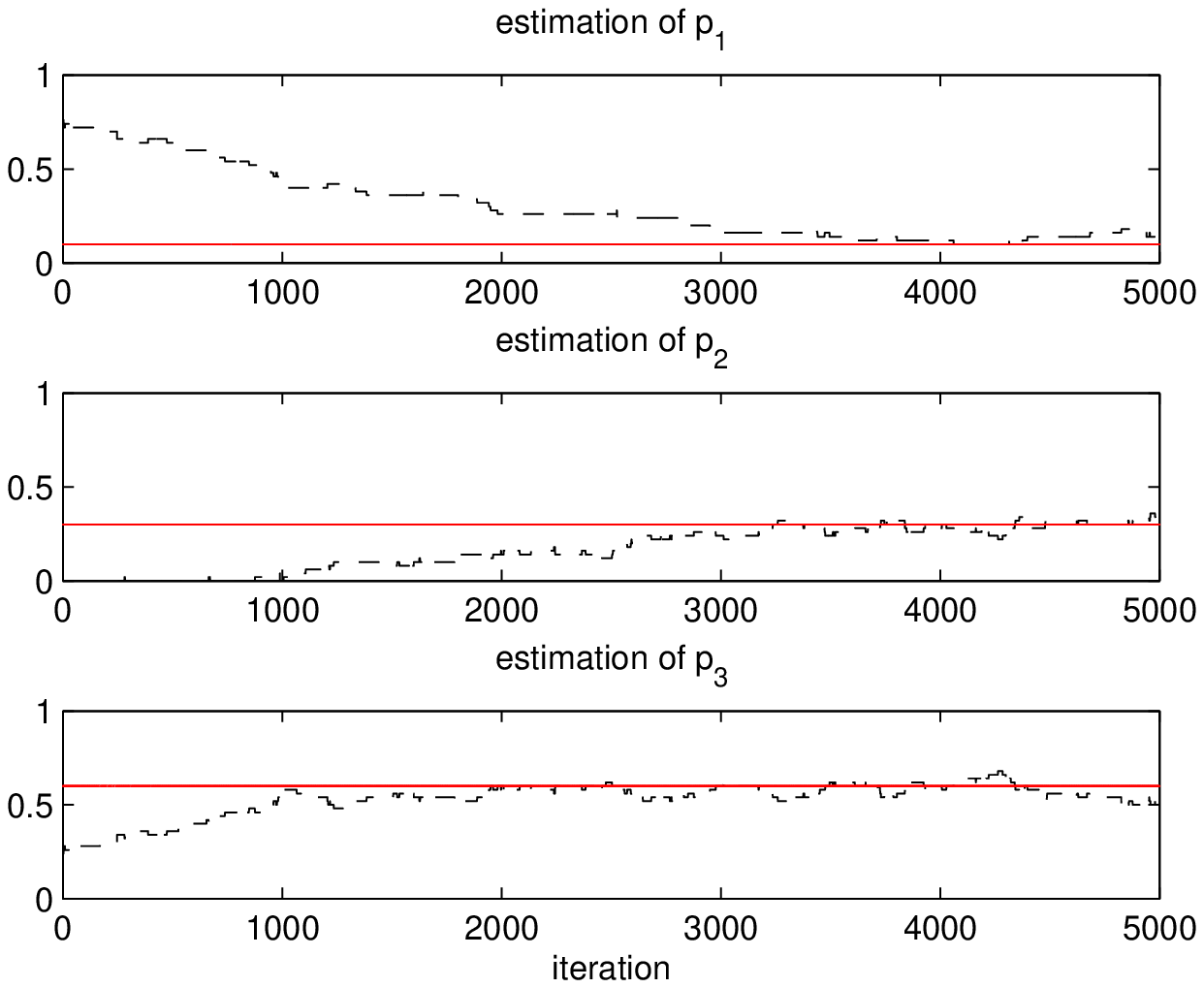}
\caption{Evolution of the proportion of particles located in the three different modes.}
\label{graphe.multi}
\end{figure}
We apply now the parallel/interacting Metropolis-Hastings sampler, see Section \ref{sec.MH}, to a case where the target distribution is multimodal:
\[
  \pi
  =
  p_1\,\NN(C_1,I) + p_2\,\NN(C_2,I) + p_3\,\NN(C_3,I)
\]
with $p_1=0.1$, $p_2=0.3$, $p_3=0.6$,  and $C_1=(-10,-10)$, $C_2=(5,0)$, $C_3=(-5,5)$. It is a mixture of 3 two-dimensional Gaussian densities.

\medskip

We describe the proposal kernel (\ref{eq.proposal}), for updating the component $X^i$, each chain $j$ propose a new candidate according to the following distribution law: 
\begin{align*}
   \piprop_{i,j}(y|X^i)
   =
   \piprop_{i,j}(y|Z^{1:i-1},X^i,X^{i+1:N})
   =
   \left\{
   \begin{array}{ll}
     \NN(X^i,{\textstyle\frac{1}{d}}\,I) \,,& \textrm{ if }i\neq j\,,
   \\[0.3em]
     \NN(X^j,I) \,,& \textrm{ if }i=j
   \end{array}
   \right.
\end{align*}
where $d\eqdef |X^i-X^j|$. 

The idea here is to explore the space with a Gaussian random walk ($i=j$) but also to allow ``jumps'' toward already explored interesting areas  
 ($i\neq j$). If $X^i$ and $X^j$ are close one the other, then ``the chain $j$ will propose a candidate far from $X^j$ and $X^i$''. If $X^i$ and $X^j$ are far one to the other, then the ``chain $j$ will propose a candidate close to $X^j$''.
 
Here $N=50$, and the initial points $\sfX^{(0),i}$, for $i=1\cdots N$, are sampled according to the uniform law on the square $[-15,10]\times [0,10]$, see Figure \ref{graphe.multi.target} (right). Figures \ref{graphe.multi.target2} clearly demonstrate the convergence of the method. In Figure \ref{graphe.multi} we present the evolution of the proportion of particles located in the neighbor of the three different modes: this also demonstrates the good behavior of the method. Note that the initial particles do not cover the mode number 2, so the algorithm is able to reach the isolated mode and to balances the particles among the modes according to the parameters $p_i$.

\subsection{An hidden Markov model}
\label{sec.hmm}

\begin{figure}[p]
\centering
\includegraphics[width=4in]{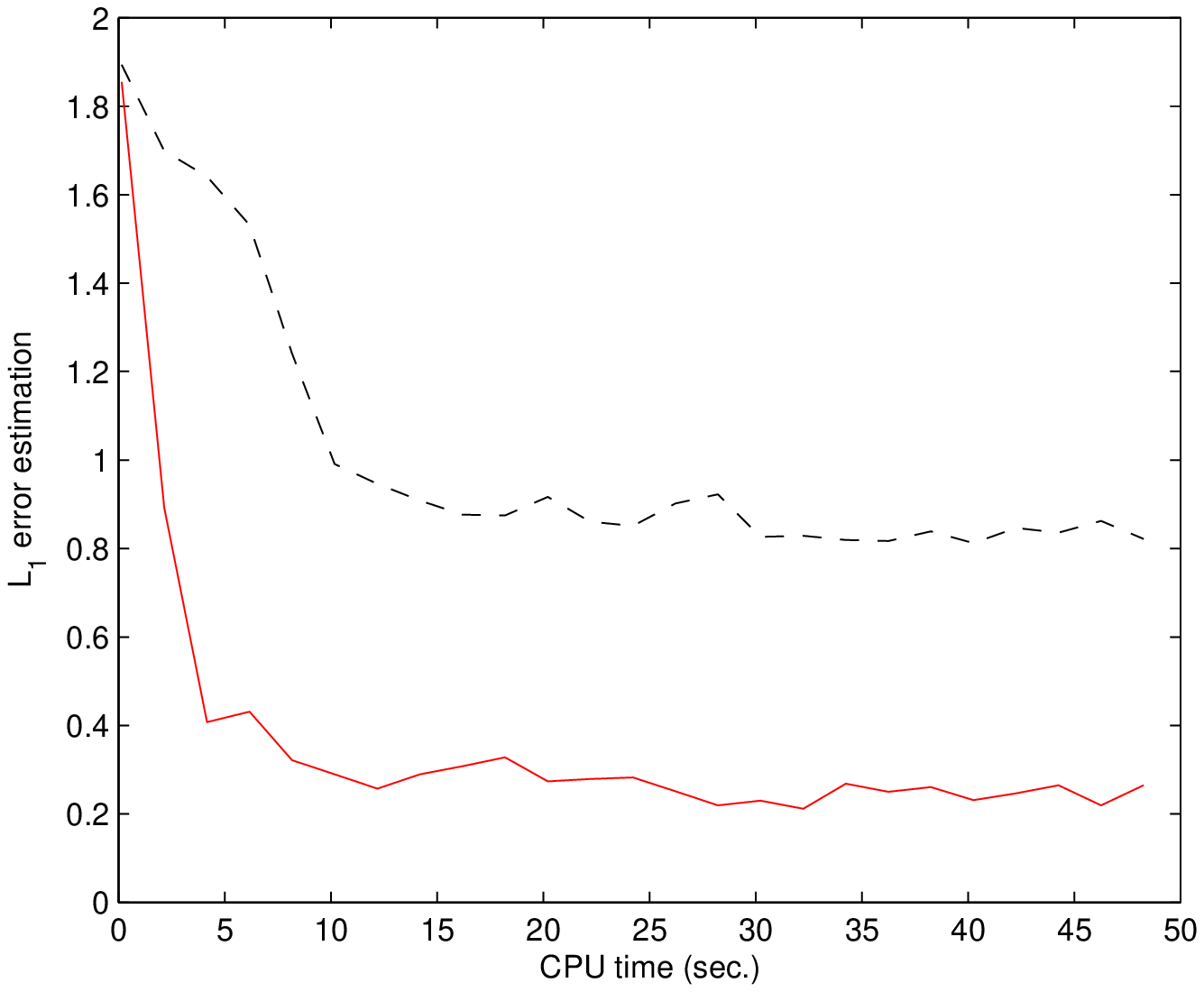}
\caption{Evolution of the indicator $\epsilon^k$, see (\ref{eq.epsi}), for the parallel/independent MwG sampler (- -), and for the parallel/interacting MH sampler (--). This evolution is depicted as a function of the CPU time and not as a function of the iteration number $k$. The residual error of about 0.22 for the second method is due to the limited size of the sample.}
\label{parallel.vs.interaction}
\includegraphics[width=4in]{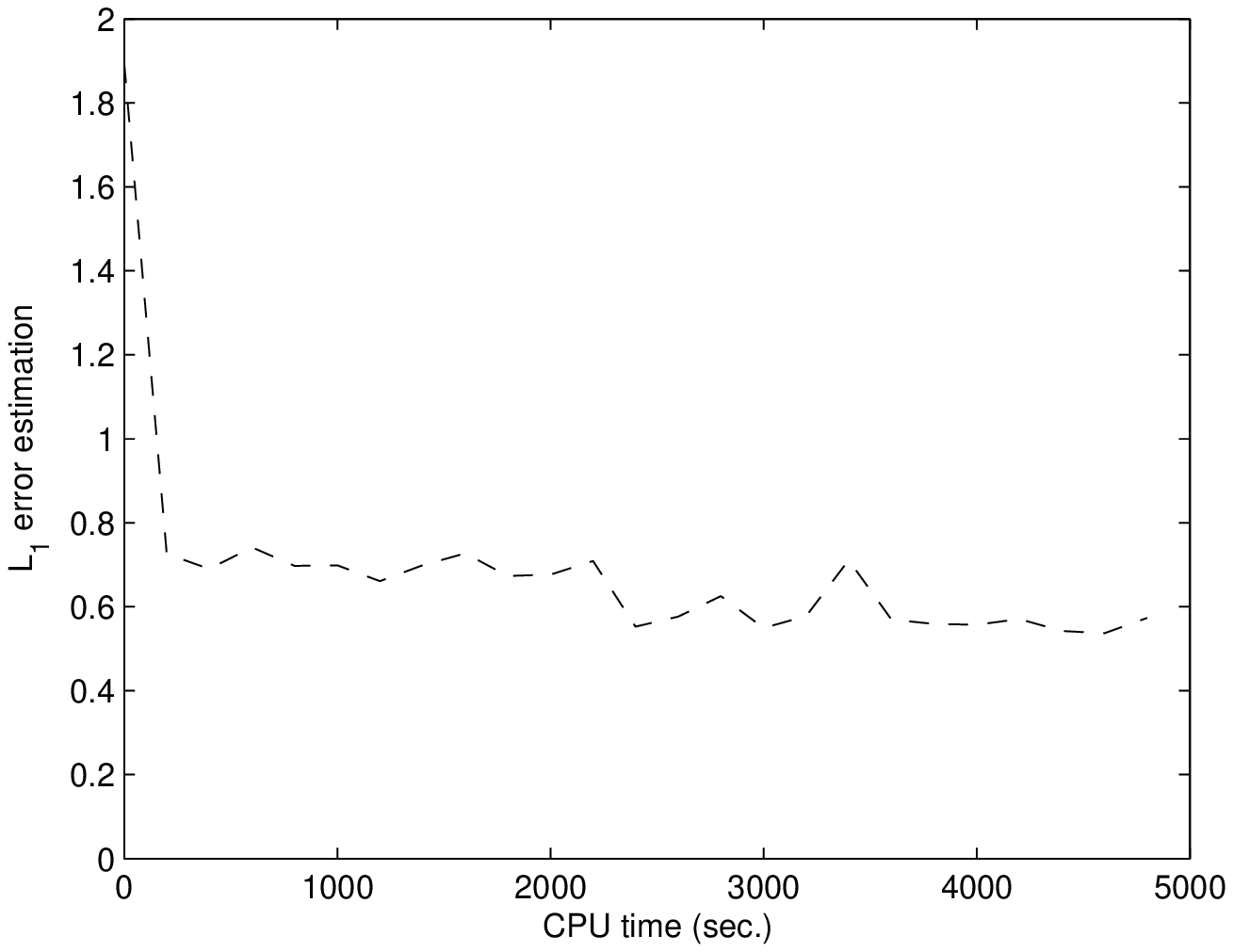}
\caption{Evolution of the indicator $\epsilon^k$, see (\ref{eq.epsi}), for the parallel/independent MwG sampler (- -). After 5000 sec. CPU time, the convergence of this method is still unsatisfactory.
}
\label{parallel.long}
\end{figure}

We apply the parallel/interacting Metropolis within Gibbs sampler, see Section \ref{sec.para.int.MwG}, to a toy problem where a good estimate $\hat{\pi}$ of the target distribution $\pi$ is available. Consider the linear Gaussian state space model:
\begin{subequations}
\label{eq.LG}
\begin{align}
\label{eq.LG.etat}
  \sfs_{\ell+1}
  &= 
  \bfa \,\sfs_\ell  +\sfw_\ell
  \,,
  &
  \sfs_1\sim \NN(\bar\sfs_1,Q_1)
  \,,
  \\
\label{eq.LG.obs}
  \sfy_\ell
  &= 
  \bfb \,\sfs_\ell +\sfv_\ell
\end{align}
\end{subequations}
for $\ell=1\cdots n$, where $\sfw_{1:n}$ and $\sfv_{1:n}$ are centered white Gaussian noises with variances $\sigma^2_\sfw$ and $\sigma^2_\sfv$. Suppose that $\bfb$ is known and $\bfa=\theta$ is unknown with a priori law $\NN(\mu_\theta,\sigma^2_\theta)$. We also suppose that $\sfw_{1:n}$, $\sfv_{1:n}$, $\sfs_1$ and $\theta$ are mutually independent. 

\medskip

The state variable is
\[
  \sfx_{1:n+1}
  \eqdef
  (\sfs_{1:n},\theta)
\]
and the target conditional density is
\[
  \pi(x_{1:n+1})
  \,\rmd x_{1:n+1}
  =
  \pi(s_{1:n},\vartheta)
  \,\rmd s_{1:n}\,\rmd\vartheta
  \eqdef
  \law(\sfs_{1:n},\theta|\sfy_{1:n}=y_{1:n})
  \,.
\]
This target law is not Gaussian, but we can perform a Gibbs sampler. Indeed the marginal conditional laws are available:
\begin{align*}
   \pi_{\sfs_\ell}(s_\ell|s_{\neg \ell},\vartheta)
   \,\rmd s_\ell
   &\eqdef
   \law(\sfs_\ell|\sfs_{\neg\ell}=s_{\neg \ell},\theta=\vartheta,\sfy_{1:n}=y_{1:n})
   =
   \NN(\mathrm{m}_\ell,\mathrm{r}^2)
   \,,
\\
  \pi_\theta(\vartheta|s_{1:n})\,\rmd \vartheta
  &\eqdef
   \law(\theta|\sfs_{1:n}=s_{1:n},\sfy_{1:n}=y_{1:n})
  =
  \NN(\mathrm{\mathrm{\tilde m},\tilde r}^2)
\end{align*}
with
\begin{align*}
  \mathrm{r}^2
  &\eqdef
  \textstyle
  \big(
    \frac{\bfb^2}{\sigma^2_\sfv}
      +\frac{1}{\sigma^2_\sfw}
      +\frac{\vartheta^2}{\sigma^2_\sfw}
  \big)^{-1}
  \,,
  &
  \mathrm{m}_\ell
  &\eqdef
  \textstyle
  \mathrm{r}^2
  \;
  \big(
    \frac{\bfb\,y_\ell}{\sigma^2_\sfv}
      +\frac{\vartheta\,s_{\ell+1}}{\sigma^2_\sfw}
      +\frac{\vartheta\,s_{\ell-1}}{\sigma^2_\sfw}
  \big)
  \,,
\\
  \mathrm{\tilde r}^2
  &\eqdef
  \textstyle
  \big(
    \frac{1}{\sigma^2_\theta}
      +\frac{\sum_{\ell=2}^n s^2_{\ell-1}}{\sigma^2_\sfw}
  \big)^{-1}
  \,,
  &
  \mathrm{\tilde m}
  &\eqdef
  \textstyle
  \mathrm{\tilde r}^2
  \;
  \big(
    \frac{\mu_\theta}{\sigma^2_\theta}
      +\frac{\sum_{\ell=2}^n s_{\ell-1}\,s_\ell}{\sigma^2_\sfw}
  \big)
  \,.
\end{align*}
We will perform three algorithms:
\begin{enumerate}
\item
$N$ parallel/interacting Metropolis within Gibbs samplers (Alg.~\ref{algo.mh.gibbs.1}),
\item
$N$ parallel/independent Metropolis within Gibbs samplers (Alg.~\ref{algo.metropolis.gibbs}),
\item
$N_{\hbox{\tiny Gibbs}}$ parallel/independent Gibbs samplers.
\end{enumerate}
Our aim is to show that making parallel samplers interact could speed up the convergence toward  the stationary distribution.

Because of its good convergence property, method \fenumiii\ is considered as a reference method. Here we perform $k=10000$ iterations of $N_{\hbox{\tiny Gibbs}}=5000$ independent Gibbs samplers. We obtain a kernel density estimate $\hat{\pi}$ of the target density based on the $N_{\hbox{\tiny Gibbs}}=5000$ final values. Let $\hat{\pi}_{\sfx_\ell}$ be the corresponding $\ell$-th marginal density. 

For methods \fenumi\ and \fenumii\ we perform $N=50$ parallel samplers. Let  $\pi^{\hbox{\tiny int},k}$ and $\pi^{\hbox{\tiny ind},k}$ be the kernel density estimates of the target density based on the final values of methods \fenumi\ and \fenumii\ respectively. Let  $\pi^{\hbox{\tiny int},k}_{\sfx_\ell}$ and $\pi^{\hbox{\tiny ind},k}_{\sfx_\ell}$ be the corresponding $\ell$-th marginal densities.

The parameter values for the simulations are $\bfa=2$, $\bfb=2$, $\sigma^2_\sfw=9$, $\sigma^2_\sfv=25$, $\sfs_1\sim \NN(4,9)$, $\theta\sim \NN(1,4)$ and $n=10$.

For each algorithm \fenumi\ and \fenumii, that is for $\pi_{\sfx_\ell}^{k}=\pi_{\sfx_\ell}^{\hbox{\tiny ind},k}$ and $\pi_{\sfx_\ell}^{\hbox{\tiny int},k}$, we compute
\[
  \epsilon_{\ell}^k
  \eqdef
  \int   
     | \pi_{\sfx_\ell}^{k}(\xi)- \hat{\pi}_{\sfx_\ell}(\xi) |
     \, \rmd \xi
  \,,\quad\ell=1\cdots n+1
  \,.
\]
Hence $\epsilon^k_\ell$ is an estimation of the $L^1$ error between the target probability distribution and its estimation provided by the algorithm used. 
To sum up the information of the $n=10$ indicators we consider their mean:
\begin{align}
\label{eq.epsi}
  \epsilon^k=\frac{1}{n+1}\sum_{\ell=1}^{n+1} \epsilon^k_\ell
  \,.
\end{align}
These estimations are based on a sample of size $N=50$ only, so they suffer from variability. This is not problematical, indeed we do not want to estimate $L^1$ errors but to diagnose the convergence toward the stationary distribution. So we use  $\epsilon^k_\ell$ as an indicator which must decrease and remain close to a small value when convergence occurs. 

To compare fairly the parallel/independent MwG algorithm and the parallel/interacted MwG algorithm, we represent on Figures \ref{parallel.vs.interaction} and \ref{parallel.long} the indicator $\epsilon^k$ for each algorithm not as a function of $k$ but as a function of the CPU time.


In Figure \ref{parallel.vs.interaction} we see that even if one iteration of algorithm \fenumi\ needs more CPU than one of \fenumii, still the first algorithm converges more rapidly than the second one. The residual error of 0.22 is due to the limited size of the sample. This error decreases to 0 as $N\uparrow\infty$. Figure \ref{parallel.long} shows the inefficiency of parallel/independent MwG on this simple model.

\section{Conclusion}
\label{conclusion}

This work showed that making parallel MCMC chains interact could improve their convergence properties.  We proved  the basic properties of the MCMC method, we did not prove that the proposed strategy speeds up the convergence. This difficult point is related to the problem of the rate of the convergence of the MCMC algorithms.

Through a simple example we saw that the Metropolis within Gibbs strategy could be a poor strategy. However this method is widely used in practice on more complex non linear models. In this situation our strategy improved the convergence properties. We also demonstrated that this approach can handle multimodal cases.

\clearpage
\appendix

\section*{Appendix: MwG algorithm}

One iteration $X\to Z$ of the Metropolis within Gibbs method consists in updating  the components $X_\ell$ successively for $\ell=1,\dots,n$, i.e.
\[
  [X_{1:n}]
  \to
  [Z_1 X_{2:n}]
  \to
  [Z_{1:2} X_{3:n}]
  \cdots
  [Z_{1:n-1} X_{n}]
  \to
  [Z_{1:n}]
  \,.
\]
Each components  $X_\ell$ is updated in two steps:
\begin{enumerate}

\item
{\it  Proposal step:}
We sample a candidate  $Y_\ell$ according to:
\[
     Y_\ell 
     \sim \piprop_\ell(\xi)\,\rmd \xi
\]

\item
{\it Selection step:}
The component $X_\ell$ could be replaced by the candidate $Y_\ell$ or stay unchanged according to a binomial sampling, the resulting value is called $Z_\ell$, i.e.:
\[
  Z_\ell
  \ot
  \left\{\begin{array}{ll}
     Y_\ell 
     & \textrm{with probability }
       \,\alpha_\ell(X_\ell,Y_\ell)  
       \,,
     \\[0.5em]
     X_\ell 
     & \textrm{with probability }
        \,1-\alpha_\ell(X_\ell,Y_\ell)
  \end{array}\right.
\]
where:
\begin{align*}
  \alpha_\ell(\xi,\xi')
  &\eqdef
  \frac{\pi_\ell(\xi')}
       {\pi_\ell(\xi )}
  \,
  \frac{\piprop_\ell(\xi )}
       {\piprop_\ell(\xi')}
  \wedge 1
\end{align*}

\end{enumerate}
The resulting algorithm is depicted in Algorithm~\ref{algo.metropolis.gibbs}.

\begin{algorithm}
\caption{\sl Metropolis within Gibbs sampler. We can go through the component indices in a random way.}
\label{algo.metropolis.gibbs}
\begin{center}
\begin{minipage}{10cm}
\hrulefill\\[-1em]
\mbox{}
\begin{algorithmic}
\STATE choose $X_{1:n}\in\R^n$ 
\FOR {$k=1,2,\dots$}
  \FOR {$\ell=1:n$}
    \STATE $Y_\ell\sim \piprop_\ell(\xi)\,\rmd \xi$ 
    \COMMENT{proposed candidate}
    \STATE $u\sim \UU[0,1]$
    \IF{$u\leq \alpha_\ell(X_\ell,Y_\ell)$} 
      \STATE $X_\ell  \ot   Y_\ell $
    \ENDIF
  \ENDFOR
\ENDFOR
\end{algorithmic}
\hrulefill
\end{minipage}
\end{center}
\end{algorithm}



\end{document}